\documentclass[10pt]{amsart}
\usepackage{amssymb}
\usepackage{amsthm}
\usepackage[longnamesfirst,numbers]{natbib}

\theoremstyle{plain}
\newtheorem{theorem}{Theorem}[section]
\newtheorem{lemma}[theorem]{Lemma}

\newtheorem{claim}[theorem]{Claim}
\newtheorem{corollary}[theorem]{Corollary}

\theoremstyle{definition}
\newtheorem{definition}{Definition}[section]

\newtheorem{example}{Example}[section]
\newtheorem{question}{Question}[section]

\theoremstyle{remark}
\newtheorem*{remark}{Remark}

\newtheorem{case}{Case}

\newcommand{\restr}{\mathrel{\mbox{\raisebox{.5mm}{$\upharpoonright $}}}}

\begin{document}
\bibliographystyle{abbrvnat}

\title{Low upper bounds of ideals}

\author{Anton\'in Ku\v{c}era}
\address{Department of Theoretical Computer Science and Mathematical Logic\\
  Faculty of Mathematics and Physics\\
  Charles University\\
  Malostransk\'{e} n\'{a}m. 25, 118 00 Praha 1\\
  Czech Republic}
\email{kucera@ksi.mff.cuni.cz}

\thanks{Ku\v{c}era was partially supported by the Research project
  of the Ministry of Education of the Czech Republic MSM0021620838}

\author{Theodore A. Slaman}
\address{Department of Mathematics\\
  The University of California, Berkeley\\
  Berkeley, CA 94720-3840 USA}
\email{slaman@math.berkeley.edu}

\thanks{Slaman was partially supported by NSF grant DMS-0501167.}

\thanks{Both authors are grateful for the support of the American
  Institute of Mathematics.  This project grew out of their
  discussions during the AIM workshop ``Effective Randomness,'' August
  7-11, 2006.  }
\thanks{We thank to an anonymous referee and to G. Barmpalias for useful comments on an earlier version of the paper}  

\keywords{Turing degree, K-trivial, low for random}
\subjclass[2000]{68Q30 (03D15)}

\begin{abstract}
  We show that there is a low $T$-upper bound for the class of
  $K$-trivial sets, namely those which are weak from the point of view
  of algorithmic randomness.  This result is a special case of a more
  general characterization of ideals in $\Delta^0_2$ $T$-degrees for
  which there is a low $T$-upper bound.
\end{abstract}

\maketitle

\section{Introduction}

\subsection{Background}

This paper is motivated by a question concerning the $K$-trivial
sets, namely those sets which are computationally weak from the point
of view of algorithmic randomness.  The collection of $K$-trivial sets
can be defined as consisting of exactly those sets for which
prefix-free Kolmogorov complexity of initial segments grows as slowly
as possible.  However, there are at least three conceptually other
ways to come to the same class, which is part of the interest in them.

Another part of the interest in this class lies in its properties
when viewed as a subideal within the $\Delta^0_2$ Turing degrees.
Here, the $K$-trivial sets induce a $\Sigma^0_3$ ideal in the
$\omega$-r.e.\ $T$-degrees which is generated by its r.e.\ members,
and the r.e.\ $K$-trivial sets induce a $\Sigma^0_3$ ideal in the
r.e.\ $T$-degrees.  This was proved by \citet{Nies:2005} and partially
also by \citet{Downey.Hirschfeldt.ea:2003} (see also
\citep{Downey.Hirschfeldt:nd} or \citep{Nies:nd*d}).  Nies
(unpublished, see \citep{Downey.Hirschfeldt:nd}) also showed that
there is a low$_2$ r.e.\ $T$-degree which is a $T$-upper bound for the
class of $K$-trivial sets.  However, \citet{Nies:2006} also proved
that there is no low r.e.\ $T$-upper bound for this class.  Since all
$K$-trivial sets are low, the latter result shows that the ideal is
nonprincipal.  Whether there is a low $T$-upper bound for the class of
$K$-trivial sets remained unresolved.  See, e.g., the list of open
questions in \citet{Miller.Nies:2006}).

This is the question which motivated this paper.  We show that there
is a low $T$-upper bound on the ideal of the $K$-trivial $T$-degrees.  The
proof applies more broadly, and we give a general characterization of
those ideals in the $\Delta^0_2$ $T$-degrees for which there is a low
$T$-upper bound.

\subsection{Notation} 

Our computability-theoretic notation generally follows
\citet{Soare:1987} and \citet{Odifreddi:1989,Odifreddi:1999}.  An
introduction to algorithmic randomness can be found in
\citet{Li.Vitanyi:1997}.  A short survey of it is also given in
\citet{Ambos-Spies.Kucera:2000} and a longer one in
\citet{Downey.Hirschfeldt.ea:2006}.  More recent progress is described
in detail in forthcoming books of \citet{Downey.Hirschfeldt:nd} and
\citet{Nies:nd*d}.

We refer to the elements of $2^{\omega}$ as sets or infinite binary
sequences.  We denote the collection of strings, i.e. finite initial
segments of sets, by $2^{<\omega}$.  The length of a string $\sigma$
is denoted by $|\sigma|$, for a set $X$, we denote the string
consisting of the first $n$ bits of $X$ by $X \restr n$.  We let
$\sigma * \tau$ denote the concatenation of $\sigma$ and $\tau$.  We
write $\sigma \preceq \tau$ to indicate that $\sigma$ is a substring
of $\tau$, and similarly $\sigma \prec \tau$ to indicate that $\sigma$
is a proper substring of $\tau$.  We further write $\sigma \prec X$ to
indicate $X\restr|\sigma| = \sigma$.  If $\sigma \in 2^{<\omega}$,
then $[\sigma]$ denotes $\{X \in 2^\omega : \sigma \prec X \}$, and
$Ext(\sigma) = \{\tau \in 2^{< \omega}: \sigma \preceq \tau \}$.

A $\Sigma^0_1$ class is a collection of sets that can be effectively
enumerated.  Such a class can be represented as $\bigcup_{\sigma \in
  W}[\sigma]$ for some (prefix-free) recursively enumerable (r.e.) set
of strings $W$.  The complements of $\Sigma^0_1$ classes are called
$\Pi^0_1$ classes. Any $\Pi^0_1$ class can be represented by the class
of all infinite paths through some recursive tree $\subseteq
2^{<\omega}$. We use also relativized versions, i.e.  $\Pi^{0,X}_1$
classes.

A string $\sigma$ is $\omega$-extendable on a tree $T$ if $\sigma
\prec X$ for some infinite path $X$ through $T$.  A string $\sigma$ is
$h$-extendable on a tree $T$ if there is a string $\tau \in T$ which
extends $\sigma$, and $|\tau|$ is equal to $|\sigma| + h$.  By
$\omega$-extendability of a string $\sigma \in 2^{<\omega}$ in a class
$\mathcal B \subseteq 2^\omega$ we mean that there is a function $f
\in \mathcal B$ which extends $\sigma$.

$[Tr]$ denotes the class of all infinite paths through a tree $Tr$.

Additionally to binary strings and trees $\subseteq 2^{<\omega}$ we
will also use both finite sequences of elements of $\{0,1\} \times
\{-1,0,1\}$ which we call $p$-strings and trees of $p$-strings
(i.e. trees $\subseteq (\{0,1\} \times \{-1,0,1\})^{<\omega}$).  We
generalize notation for binary strings to $p$-strings in an obvious
way.

$K$ denotes prefix-free Kolmogorov complexity.  We assume
Martin-L\"{o}f's definition of $1$-randomness as well as its
relativization to an oracle.

Let $f$ be a total function, i.e. $f \in \omega^\omega$, and let $g$ be a 
(possibly) partial function, $g : \omega  \mapsto \omega$. 
We say that $f$ eventually dominates $g$ if there is an $n$ 
such that for all $x > n$, \ $g(x)\downarrow$ \ implies \ $f(x) \geq g(x)$.

\textit{Convention.}  If $\mathcal C$ is a class of sets, $T$-degrees
of which form an ideal, we often speak of an ideal $\mathcal C$ by
which we mean an ideal of $T$-degrees of members of $\mathcal C$.

\section{Preliminaries}

We begin with basic definitions with which one can formulate the
various characterizations of $K$-triviality.  However, these
characterizations differ in applicability as will be mentioned later.

\begin{definition}
  \begin{enumerate}
  \item ${\mathcal K}$ denotes the class of $K$-trivial sets, i.e. the
    class of sets $A$ for which there is a constant $c$ such that for
    all $n,$ $K(A \restr n) \leq K(0^n) + c$.
  \item $\mathcal L$ denotes the class of sets which are low for
    $1$-randomness, i.e. sets $A$ such that every $1$-random set is
    also $1$-random relative to $A$.
  \item $\mathcal M$ denotes the class of sets that are low for $K$,
    i.e. the class of sets $A$ for which there is a constant $c$ such
    that for all $\sigma$, $K(\sigma) \leq K^A(\sigma) + c.$
  \item A set $A$ is a basis for $1$-randomness if there is a $Z$ such
    that $A \leq_T Z$ and $Z$ is $1$-random relative to $A$. The
    collection of such sets is denoted by $\mathcal {BR}$.
  \end{enumerate}
\end{definition}

\citet{Nies:2005} proved that $\mathcal L$ = $\mathcal M$, Hirschfeldt
and Nies, see \citep{Nies:2005}, proved that $\mathcal K$ = $\mathcal
M$, and \citet{Hirschfeldt.Nies.ea:nd} proved that $\mathcal {BR} =
\mathcal K$.  Thus, all these four classes are equal and we have,
remarkably, four different characterizations of the same class.  

\citet{Chaitin:1977} proved that if a set is $K$-trivial then it is
$\Delta^0_2$.  By a result of \citet{Kucera:1993}, low for $1$-random
sets are $GL_1$ and, thus, every $K$-trivial set is low. The lowness
of these sets also follows from some recent results on this class of
sets, see \citep{Nies:2005} or \citep{Downey.Hirschfeldt.ea:2006}.

An interesting result shows that there is an effective listing of all
$K$-trivial sets in the following way.

\begin{theorem}[\citet*{Downey.Hirschfeldt.ea:2003}]
  There is an effective sequence $\{B_e,d_e\}_e$ of all the r.e.\
  $K$-trivial sets and of constants such that each $B_e$ is
  $K$-trivial via $d_e$.
\end{theorem} 

\citet{Nies:2005} proved that the class of $K$-trivial sets is closed
downwards under $T$-reducibility.  Of course, the downward closure is
immediate for the equality $\mathcal K=BR$, but that equality is a
deeper fact.  In the same paper, Nies also showed that for any
$K$-trivial set $A$ there is an r.e.\ $K$-trivial set $B$ such that $A
\leq_{tt} B$.  \citet{Downey.Hirschfeldt.ea:2003} proved that the
class of $K$-trivial sets is closed under join.  Finally,
\citet{Nies:2005} showed that the $K$-trivial sets induce a
$\Sigma^0_3$ ideal in the $\omega$-r.e.\ $T$-degrees which is
generated by its r.e.\ members, and r.e.\ $K$-trivial sets induce a
$\Sigma^0_3$ ideal in the r.e.\ $T$-degrees.

On the other hand, the ideal is nonprincipal as the following theorem
shows, since all $K$-trivials are low.

\begin{theorem} [\citet{Nies:2006}]\label{2.2}
  \begin{itemize}
  \item For each low r.e.\ set $B$, there is an r.e.\ $K$-trivial set
    $A$ such that $A \nleq_T B$.
  \item For any effective listing $\{B_e,z_e\}_e$ of low r.e.\ sets and
    of their low indices there is an r.e.\ $K$-trivial set $A$ such
    that $A \nleq_T B_e$ for all $e$.
  \end{itemize}
\end{theorem}

\begin{remark}
  The proof uses a technique, known as Robinson low guessing method
  (introduced for low r.e.\ sets by \citet{Robinson:1971}) which is
  compatible for low r.e.\ sets with a technique {\it do what is
    cheap}.  Here {\it cheap} is defined in terms of a cost function
  (see, e.g. \citep{Downey.Hirschfeldt:nd},
  \citep{Downey.Hirschfeldt.ea:2006}, \citep{Nies:nd*d}).
  Alternatively, one could construct in the above Theorem a set which
  is low for $1$-randomness instead of a $K$-trivial set and define
  {\it cheap} as having a small measure (see,
  e.g. \citet{Kucera.Terwijn:1999} or
  \citet{Downey.Hirschfeldt.ea:2006}).  However, the Robinson low
  guessing technique does not seem to generalize from r.e.\ sets to
  $\Delta^0_2$ sets in a way which is compatible with the heuristic
  {\it do what is cheap}. In fact, it does not as Theorem~\ref{3.2}
  below shows.
\end{remark}

An immediate corollary of Theorem~\ref{2.2} is that no low r.e.\ set
can be a $T$-upper bound for the class $\mathcal K$.

Further, Theorem~\ref{2.2} was used by \citet{Nies:2006, Nies:2005}
and by \citet{Downey.Hirschfeldt.ea:2003} to show that four different
characterizations of the same class, i.e. characterizations yielding
$\mathcal K, \mathcal L, \mathcal M, \mathcal {BR}$ respectively, are
not equally uniform.  Especially, the uniformity in the
characterization by $K$-triviality sets is weaker than the
characterization by low for $K$.  Unlike the constants by which a set
is $K$-trivial, constants by which a set is low for $K$ can be
uniformly transformed into indices by which that set is low.

\begin{theorem}[\citet{Nies:2006}, \citet{Nies:2005},
  \citet*{Downey.Hirschfeldt.ea:2003}]\label{2.3}\hfill\mbox{}
  \begin{itemize} \item There is no effective sequence
      $\{B_e,c_e\}_e$ of all the r.e.\ low for $K$ sets with
      appropriate constants.
    \item There is no effective way to
      obtain from a pair $(B,d)$, where $B$ is an r.e.\ set that is
      $K$-trivial via $d$, a constant $c$ such that $B$ is low for $K$
      via $c$.
    \item There is no effective listing of all the r.e.\
      $K$-trivial sets together with their low
      indices.
    \end{itemize}
  \end{theorem}

  There are several results
  on $\Sigma^0_3$ ideals of r.e.\ sets.  One of the first such results
  is the following.

  \begin{theorem}[\citet{Yates:1969}] For any r.e
    set $A <_T \emptyset'$ the following conditions are
    equivalent.
    \begin{enumerate}
    \item $A'' \equiv_T \emptyset''$.
    \item $\{x : W_x \leq_T A \}$ is a $\Sigma^0_3$ set.
    \item the class $\{W_x : W_x \leq_T A \}$ is uniformly
      r.e.
    \end{enumerate}
  \end{theorem}

  On the other hand Nies (unpublished, see
  \citep{Downey.Hirschfeldt:nd}) proved that for any proper
  $\Sigma^0_3$ ideal of r.e.\ sets there is a low$_2$ r.e.\ $T$-degree
  which is a $T$-upper bound for this ideal.  Thus, together with the
  above theorem by Yates we have a characterization of ideals of r.e.\
  sets for which there is a low$_2$ r.e.\ $T$-upper bound.  An ideal
  has a low$_2$ $T$-upper bound if and only if it is a subideal of a
  proper $\Sigma^0_3$ ideal.

  A characterization of ideals of r.e.\ sets (or ideals generated by
  their r.e.\ members) for which there is a low $T$-upper bound, not necessarily
  r.e., was open.  We substantially use properties of $\{0,1\}$-valued
  DNR functions (and their relativizations) in our construction of low
  $T$-upper bounds for ideals.  We give a short review of basic properties
  of such functions here.

  \begin{definition} Let $\mathcal {PA}(B)$ denote the class of all
    $\{0,1\}$-valued $B$-DNR functions, i.e.  the class of functions
    $f \in 2^\omega$ such that $f(x) \neq \Phi_x(B)(x)$ for all $x$ .
    If $B$ is $\emptyset$ we simply speak of $\mathcal {PA}$.
  \end{definition}

  \begin{definition}[\citet{Simpson:1977}]
    Write ${\bf b} << {\bf a}$ to mean that every infinite tree $T
    \subseteq 2^{<\omega}$ of $T$-degree $\leq {\bf b}$ has an
    infinite path of $T$-degree $\leq {\bf a}$.
  \end{definition}

  \begin{theorem}[\citet{Scott:1962}, Solovay (unpublished),
    see \citep{Simpson:1977}]\label{2.5} 
    The following conditions are equivalent:
    \begin{enumerate}
    \item
      ${\bf a}$ is a $T$-degree of a $\{0,1\}$-valued DNR
      function.  \item ${\bf a} >> {\bf 0}$.
    \item ${\bf a}$ is a $T$-degree of a complete extension of Peano
      arithmetic.
    \item ${\bf a}$ is a $T$-degree of a set separating some effectively
      inseparable pair of r.e.\ sets.
    \end{enumerate}
  \end{theorem}

  \begin{remark}
    By the implication from (1) to (2) in Theorem~\ref{2.5}, $\mathcal
    {PA}$ is a ``universal'' $\Pi_1^0$ class.  $\{0,1\}$-valued DNR
    functions are also called PA sets and $T$-degrees $>> {\bf 0}$ are
    also called PA degrees.  Analogously, the class $\mathcal {PA}(B)$
    is a ``universal'' $\Pi_1^{0,B}$ class.  \citet{Simpson:1977}
    proved that the partial ordering $<<$ is dense and $\bf a << \bf
    b$ implies $\bf a < \bf b$.
  \end{remark}
  
\begin{definition} Let
    $M$ be an infinite set and $\{m_i\}_i$ be an increasing list of
    all members of $M$.  For $f \in 2^\omega$ by $Restr(f,M)$ we
    denote a function $g$ defined for all~$i$ by $g(i) = f(m_i)$.
    Similarly, if $\mathcal B \subseteq 2^\omega$ then by
    $Restr(\mathcal B,M)$ we denote a class of functions $\{g: g =
    Restr(f,M) \ \& \ f \in \mathcal B\}$. Further, if $\sigma \in
    2^{< \omega}$, then by $Restr(\sigma,M)$ we denote a string $\tau$
    defined by $\tau(i) = \sigma(m_i)$ for all $i$ such that $|m_i| <
    |\sigma|$.
  \end{definition}

  \begin{lemma}\label{2.6}
    \begin{enumerate}
    \item For every $\Pi^0_1$ class $\mathcal B$ which is a subclass
      of $\mathcal {PA}$ there is an infinite recursive set $M$ such
      that if $\mathcal B$ is nonempty then $Restr(\mathcal B,M) =
      2^\omega$, i.e. for every function $g \in 2^\omega$ there is a
      function $f \in \mathcal B$ such that $Restr(f,M) =
      g$. Moreover, an index of $M$ can be found recursively from an
      index of $\mathcal B$.
    \item For every $\Pi^{0,B}_1$ class $\mathcal B$ which is a
      subclass of $\mathcal {PA}(B)$ there is an infinite recursive
      set $M$ such that if $\mathcal B$ is nonempty then
      $Restr(\mathcal B,M) = 2^\omega$, where an index of $M$ can be
      found uniformly-recursively from an index of $\mathcal B$,
      i.e.\ in a uniform way which does not depend on an oracle $B$.
    \end{enumerate}
  \end{lemma}


\begin{proof} Lemma~\ref{2.6} is an application of G\"{o}del's
 incompleteness phenomenon in the context of $\Pi^0_1$ classes of
 $\{0,1\}$-valued $B$-DNR functions.  Under a slight modification
 it was proved by \citet{Kucera:1989}. For the convenience of the
 reader we first sketch the proof of the original version from \citep{Kucera:1989},
 which is slightly technically easier. Then we show how to modify it to get the desired result stated in Lemma~\ref{2.6}.
 
 \noindent{\it  1. The original version from \citep{Kucera:1989}.} 

 Let $\mathcal {G}_0$ denote the class of all  $\{0,1\}$-valued GNR functions, i.e. functions
 $f \in 2^\omega$ such that $f(<x,y>) \neq \varphi_x(y)$ for all $x,y$.
 Theorem 2 from \citep{Kucera:1989} states the following.\\
 For every nonempty $\Pi^0_1$ class $\mathcal B$, $\mathcal B \subseteq \mathcal {G}_0$, there is $x_0$ such that
 for every set $C$ there is a function $g \in \mathcal B$ such that 
 $g(<x_0,y>) = C(y)$ for all $y$ (and such $x_0$ can be found effectively from an index of $\mathcal B$).\\
 Proof of this theorem.\\
 Suppose $\mathcal B$ is a $\Pi^0_1$ subclass  of $\mathcal {G}_0$.
 Observe that $\mathcal B \cap \{f \in 2^\omega : f(<x,j>) = \sigma(j), \ j < |\sigma| \} = \emptyset$
 is a $\Sigma^0_1$ condition (with variables $x$ and $\sigma \in 2^{<\omega}$). Thus, there is a total recursive function $\beta$ such that 
 $\varphi_{\beta(x)}(y)$ is defined for all $y$ if and only if there is a string $\sigma \in 2^{<\omega}$ 
 for which 
 $\mathcal B \cap \{f \in 2^\omega : f(<x,j>) = \sigma(j), \ j < |\sigma| \} = \emptyset $, 
 and if there is such a string, then the first such found under a standard search, say $\sigma_0$, is used to define
 $\varphi_{\beta(x)}(y) = 1 - \sigma_0(y)$ for all $y < |\sigma_0|$ and $\varphi_{\beta(x)}(y) = 0$ for all $y \geq |\sigma_0|$.
 By the recursion theorem there is $x_0$ such that $\varphi_{x_0} = \varphi_{\beta(x_0)}$.\\
 It is easy to prove that if $\mathcal B$ is nonempty, then both 
 $\mathcal B \cap \{f \in 2^\omega : f(<x_0,j>) = \sigma(j), \ j < |\sigma| \} \neq \emptyset $ for any   $\sigma \in 2^{<\omega}$
 and $\varphi_{x_0}(y)$ is not defined for any $y$. Really, it there were such a string $\sigma$ 
 for which $\mathcal B \cap \{f \in 2^\omega : f(<x_0,j>) = \sigma(j), \ j < |\sigma| \} = \emptyset $ and
 $\sigma_0$ were the first such found (as mentioned above) then every function $g \in \mathcal {G}_0$ would satisfy 
 $g(<x_0,j>) = \sigma_0(j)$ for $j < |\sigma_0|$, which would immediately imply $\mathcal B = \emptyset$,
 a contradiction. Thus, for every string $\sigma \in 2^{< \omega}$ we have
 $\mathcal B \cap \{f \in 2^\omega : f(<x_0,j>) = \sigma(j), \ j < |\sigma| \} \neq \emptyset $.
 Using the compactness of the Cantor space $2^\omega$ we have the required property, i.e. 
 $\mathcal B \cap \{f \in 2^\omega : f(<x_0,y>) = C(y), \ y \in \omega  \} \neq \emptyset$  for any set $C$.
 
\noindent{\it 2. Proof of the version stated in Lemma~\ref{2.6}, part 1.}

  We easily modify the above proof by an additional use of the s-m-n theorem.
  Suppose a $\Pi^0_1$ class $\mathcal B$ is a subclass of $\mathcal {PA}$.\\
  Let $\alpha$ be a partial recursive function such that 
  $\alpha(z,y,w)$ is defined if and only if
  there is a string $\sigma \in 2^{< \omega}$ such that 
  $\varphi_z(j)$ is defined for all $j < |\sigma|$ and   
  $\mathcal B \cap \{f \in 2^\omega : f(\varphi_z(j)) = \sigma(j), \ j < |\sigma| \} = \emptyset $, 
  and if such string exists then the first such found under a standard search, say $\sigma_0$, is used to define
  $\alpha(z,y,w) = 1 - \sigma_0(y)$ for all $y < |\sigma_0|$ and $\alpha(z,y,w) = 0$ for $y \geq |\sigma_0|$.
  By the s-m-n theorem there are total recursive functions $\gamma$ and $\delta$ (even primitive recursive) 
  such that 
  $\alpha(z,y,w) = \varphi_{\gamma(z,y)}(w)$  and  $\gamma(z,y) = \varphi_{\delta(z)}(y)$, for all $z,y,w$.\\
  By the recursion theorem there is $z_0$  such that $\varphi_{z_0}(y) = \varphi_{\delta(z_0)}(y) = \gamma(z_0,y)$ for all $y$.
  By properties of s-m-n functions we may assume that the function  $\gamma$ is increasing in both variables.
  Let $h$ denote $\varphi_{z_0}$. Then $h$ is a total increasing recursive function,
  and $\alpha(z_0,y,w) = \varphi_{h(y)}(w)$, for all $y,w$.
  A straightforward modification of the above proof of the original version shows that
  if $\mathcal B$  is nonempty then 
  $\mathcal B \cap  \{f \in 2^\omega : f(h(j)) = \sigma(j), \ j < |\sigma| \}\neq \emptyset$ for every string $\sigma \in 2^{< \omega}$.
  To  verify that in details, suppose for a contradiction that there is a string $\sigma \in 2^{< \omega}$ 
  for which $\mathcal B \cap \{f \in 2^\omega : f(h(j)) = \sigma(j), \ j < |\sigma| \} = \emptyset $
  and let $\sigma_0$ denote the first such string found by a search realized to define values
  of the function $\alpha$, namely values $\alpha(z_0,y,w)$.
  Then $\alpha(z_0,y,w)$ is defined and 
  $\alpha(z_0,y,w) = \varphi_{h(y)}(w) = \varphi_{h(y)}(h(y))$ for all $y,w$, and, further, 
  $\varphi_{h(y)}(h(y)) = 1 - \sigma_0(y)$ for all $y < |\sigma_0|$.
  Then, by the definition of the class $\mathcal {PA}$, i.e. the class of $\{0,1\}$-valued DNR functions ,
  every function $f \in \mathcal {PA}$ has to satisfy  $f(h(y)) \neq \varphi_{h(y)}(h(y))$ 
  and, thus, $f(h(y)) = \sigma_0(y)$ for all $y < |\sigma_0|$. This immediately yields that
  $\mathcal B = \mathcal B \cap \{f \in 2^\omega : f(h(j)) = \sigma_0(j), \ j < |\sigma_0| \} $ and, therefore, 
  this proves that the class $\mathcal B$ is empty, a contradiction.\\
  Thus, for every string $\sigma \in 2^{< \omega}$ we have
  $\mathcal B \cap \{f \in 2^\omega : f(h(j)) = \sigma(j), \ j < |\sigma| \} \neq \emptyset $.
  Using the compactness of the Cantor space $2^\omega$ we have the required property, i.e. 
  $\mathcal B \cap \{f \in 2^\omega : f(h(y)) = C(y), \ y \in \omega  \} \neq \emptyset$  for any set $C$.
  Let $M$ denote the recursive set which is the range of $h$. Obviously $M$ has all required properties.

\noindent{\it 3. Part 2 of Lemma~\ref{2.6}.}

  It is just a relativized version of part 1 and it is proved analogously.

\end{proof}

  By Lemma~\ref{2.6}, we can code arbitrary sets into members of
$\Pi^{0,B}_1$ subclasses of $\mathcal {PA}(B)$. We illustrate it for
an unrelativized case.  Suppose that $\mathcal B$ is a nonempty
$\Pi^0_1$ subclass of $\mathcal {PA}$ and $M$ an infinite recursive
set such that $Restr(\mathcal B,M) = 2^\omega$. Let a set $C$ be
given.  If we take a class $\mathcal E = \{f : f \in \mathcal B \ \& \
Restr(f,M) = C \}$, then by our assumption $\mathcal E$ is nonempty.
It is a $\Pi^{0,C}_1$ class and obviously any member $B$ of $\mathcal
E$ is $T$-above $C$. In a more general way, we may nest into a
$\Pi^0_1$ subclass of $\mathcal {PA}$ not only a singleton $\{C\}$ as
above, but even a given $\Pi^{0,C}_1$ class. Namely, with the above
assumptions if $\mathcal C$ is a nonempty $\Pi^{0,C}_1$ class then a
class $\mathcal E = \{f : f \in \mathcal B \ \& \ Restr(f,M) \in
\mathcal C \}$ is nonempty, it is a $\Pi^{0,C}_1$ class and obviously
any member of $\mathcal E$ is $T$-above some member of $\mathcal C$.
A relativization of these tricks to $\Pi^{0,B}_1$ classes which are
subclasses of $\mathcal {PA}(B)$ is straightforward.  When combined
with Low Basis Theorem of \citet{Jockusch.Soare:1972} we easily get
the following.

\begin{example}\label{2.1}
  For any low set $A$ there is a low PA set $B$ such that $A \leq_T
  B$. We may even require that additionally $\bf a << \bf b$ where
  $\bf a, \bf b$ are $T$-degrees of sets $A,B$ respectively.  A
  relativization of this fact to an oracle is straightforward.
\end{example}


\section{Constructing low upper bounds for ideals}

We show that there is a low $T$-upper bound for the class of
$K$-trivial sets and give also a more general result about low
$T$-upper bounds for ideals in the $\Delta^0_2$ $T$-degrees.

\begin{theorem}\label{3.1}
  Let $\mathcal C$ be a $\Sigma^0_3$ ideal of r.e.\ sets. Then the
  following conditions are equivalent.
  \begin{enumerate}
  \item There is a function $F$ recursive in $\emptyset'$ which
    eventually dominates all partial functions recursive in any member of
    $\mathcal C$.
  \item There is a low $T$-upper bound for $\mathcal C$.
  \end{enumerate}
\end{theorem}

Theorem ~\ref{3.1} follows from the next more general result.

\begin{theorem}\label{3.2}
  Let $\mathcal C$ be an ideal in $\Delta^0_2$ $T$-degrees.
  The following conditions (1) and (2) are equivalent.
  \begin{enumerate}
  \item
    \begin{enumerate}
    \item $\mathcal C$ is contained in an ideal $\mathcal A$ which is 
      generated by a sequence of sets $\{A_n\}_n$ such that the
      sequence  is uniformly recursive in $\emptyset'$ and
    \item there is a function $F$ recursive in $\emptyset'$ which
      eventually dominates any partial function recursive in any set
      with $T$-degree in $\mathcal A$.
    \end{enumerate}
  \item There is a low $T$-upper bound for $\mathcal C$.
  \end{enumerate}
\end{theorem}

\begin{remark} 
  We may equivalently require that a low $T$-upper bound (mentioned in
  Theorems~\ref{3.1} and~\ref{3.2}) is PA since, as we saw, every low
  set has a low PA set $T$-above it. Thus, $T$-upper bounds which are
  PA are the most general case in this characterization.
\end{remark}

\begin{corollary} 
  There is a low set which is a $T$-upper bound for the class $\mathcal K$,
  i.e. for the ideal of $K$-trivial sets.
\end{corollary}

\begin{proof} 
  As we already mentioned (see \citep{Nies:2005}) the class of r.e.\
  $K$-trivial sets induces a $\Sigma^0_3$ ideal in the r.e.\
  $T$-degrees, and the ideal $\mathcal K$ is induced by its r.e.\
  members.   \citet{Kucera.Terwijn:1999} proved that there is a
  function $F$ recursive in $\emptyset'$ which eventually dominates all partial
  functions recursive in any set which is low for $1$-randomness.
  Since $\mathcal L = \mathcal K$, the corollary follows.
\end{proof}

\begin{remark}
  We explain the main obstacles of proving Theorem ~\ref{3.2}.

  The implication from (2) to (1) is direct.  If $L$ is a low set,
  then $\emptyset'$ can compute the function $f:n\mapsto m$, where $m$ is the
  strict supremum of the set
  \[\{\{e\}^L(n): e \leq n \ \& \ \{e\}^L(n)\mbox{ converges}\}.
  \]This function eventually dominates every function recursive in
  $L$.  Similarly, $\emptyset'$ can uniformly-recursively compute a
  sequence of sets consisting of exactly those sets which are
  recursive in $L$.  For example, take the sequence $X_e$ such that
  $X_e=\{n:\{e\}^L(n)=1 \ \& \ \{e\}^L(j)\mbox{ converges for all } j
  \leq n \}$.

  The implication from (1) to (2) is more subtle.  Assume that
  $\mathcal C$ is generated by the uniformly recursive in $\emptyset'$
  sequence of sets $\{A_n\}_n$ and there is a function $F$ recursive
  in $\emptyset'$ which eventually dominates any partial function
  recursive in any set with $T$-degree in $\mathcal C$.  (We identify
  $\mathcal C$ with the ideal $\mathcal A$ described in (1).)

  We want to construct recursively in $\emptyset'$ a low set $A$ for
  which $A \geq_T A_n$ for all $n$. The obstacle is that we do not
  have low indices of sets $A_n$, but we have to effectively decide
  facts about $A'$.  A solution consists in using relativized
  $\Pi^0_1$ classes.  Once we commit ourselves to $A$'s satisfying a
  $\Pi^0_1$ sentence, we must ensure that sentence to be true in the
  limit.  In other words, our commitment is that $A$ should belong to
  the $\Pi^0_1$ class of reals for which the $\Pi^0_1$ sentence is
  true.  Our next problem is coding the given sets $A_n$ into the
  members of (relativized) $\Pi^0_1$ classes to which we have
  committed ourselves.  Here a solution consists in using "rich"
  $\Pi^0_1$ classes or relativized $\Pi^0_1$ classes, like $\mathcal
  {PA}$ or $\mathcal {PA}(A_n)$.  Finally, we come to the technical
  finesse in the construction.  We use the given function $F$ recursive
  in $\emptyset'$ which eventually dominates every partial function
  recursive in any $A_n$ to replace missing low indices of $A_n$.
  More precisely, we replace questions about $\omega$-extendability of
  a string on an $A_n$-recursive tree by questions about its
  finite-extendability where the depth to which extendability is
  required is computed by $F$.  By a finite injury construction, we
  can guarantee that eventually the answers to our questions about
  appropriate finite-extendability (of a string on a given
  $A_n$-recursive tree) computed by $F$ are, in fact, correct answers
  to questions about $\omega$-extendability.


  Before presenting the whole (global) construction of the desired set $A$
  we first explain the main idea on a simple case, where we deal with just
  one given set $A_n$ individually (in an isolated way).
  Later we describe how to combine all individual cases together to provide
  one global construction.

  We give a simpler version for the case of $K$-trivial sets first and
  then give a more general version which is needed to prove
  Theorem~\ref{3.2}.
\end{remark}

\begin{lemma}\label{3.4}
  There is a recursive procedure which given an index of an r.e.\
  $K$-trivial set $A$ produces a low set $A^*$ and the lowness index
  for $A^*$ such that $A \leq_T A^*$.  That is, there are recursive
  functions $f, g$ such that if $W_e$ is $K$-trivial then
  $\Phi_{f(e)}(\emptyset')$ is a low set, $g(e)$ is its lowness index
  and $W_e \leq_T \Phi_{f(e)}(\emptyset')$.
\end{lemma}

\begin{remark} 
  We do not claim that $A \leq_T A^*$ uniformly in an index
  of $A$. In fact, this uniformity would immediately provide
  uniformity for low indices of all r.e.\ $K$-trivial sets.  But this
  contradicts the result of \citet{Nies:2006} described in
  Theorem~\ref{2.2} or, equivalently, it contradicts
  Theorem~\ref{2.3}.  Also, by Theorem~\ref{2.2}, 
  the sets $A^*$ mentioned in Lemma~\ref{3.4}
  cannot be obtained uniformly as r.e.\ sets (i.e. presented by their
  r.e.\ indices).
\end{remark}

The previous lemma follows immediately from the more general
result stated in Lemma~\ref{3.5}.  

\begin{lemma}\label{3.5}
  Given a function $F$ recursive in $\emptyset'$, there is a uniform
  way to obtain from a $\emptyset'$-index of a set $A$ with the
  property that any partial function recursive in $A$ is eventually dominated by
  $F$ both a low set $A^*$ and an index of lowness of $A^*$ such that
  $A \leq_T A^*$.  That is to say that there are recursive functions
  $f, g$ such that if $\Phi_e(\emptyset')$ is total and equal to some
  set $A$ so that any partial function recursive in $A$ is eventually dominated
  by $F$ then $\Phi_{f(e)}(\emptyset')$ is a low set, $g(e)$ is its
  lowness index and $A \leq_T \Phi_{f(e)}(\emptyset')$.
\end{lemma}

\begin{proof}
  Lemma~\ref{3.5} is the heart of the matter and its proof is the most
  technically involved section in our analysis.
  
  The idea behind the proof is to combine forcing with $\Pi^0_1$
  classes (as in the \citet{Jockusch.Soare:1972} Low Basis Theorem)
  with coding sets into members of nonempty $\Pi^0_1$ subclasses of
  the class $\mathcal {PA}$.  The given function $F$ is used to
  approximate the answers to $A'$-questions.  If $A$ satisfies the
  given assumptions, our method will guarantee that the approximation
  will be correct from some point on.  For the reader who is steeped
  in the priority methods of the recursively enumerable Turing
  degrees, our construction is the implementation of a
  $\Sigma^0_2$-strategy (like those in the Sacks Splitting Theorem) in
  which each action by the strategy restricts the construction to a
  yet smaller $\Pi^0_1$-class.  Note, that $(A^*)'$ has to be
  uniformly recursive in $\emptyset'$. Thus, our
  $\emptyset'$-construction cannot change any decision about
  $(A^*)'(x)$ that it has already made.

  We now describe our construction for one given $A$ . So, let $A$ be
  recursive in $\emptyset'$, let $F$ be recursive in $\emptyset'$, and assume
  that any partial function recursive in $A$ is eventually dominated by $F$.

  Let $PR$ denote $\{0,1\} \times \{-1,0,1\}$ and let $PR^{< \omega}$
  denote the set of all finite sequences of elements of $PR$. We call
  such sequences as $p$-strings. We use standard notation associated
  with binary strings also for $p$-strings in an obvious way. Any
  $p$-string $\rho$ may be viewed as a pair $(\sigma,\alpha)$ of a
  binary string $\sigma \in 2^{< \omega}$ and a sequence $\alpha \in
  \{-1,0,1\}^{< \omega}$ both of lengths equal to $|\rho|$, for which
  $\rho(j) = (\sigma(j),\alpha(j))$ for all $j < |\rho|$.

  Further, for any finite sequence $\beta$ from $\{-1,0,1\}^{<
    \omega}$ let $\beta^{\geq 0}$ denote a binary string $\alpha$
  which arises by deleting all $(-1)$'s from $\beta$. Similarly, for
  any infinite sequence $X$ from $\{-1,0,1\}^{\omega}$, $X^{\geq 0}$
  denote (finite or infinite) binary sequence arising by deleting all
  $(-1)$'s from $X$.

  Working recursively in $\emptyset'$, we construct an infinite
  perfect tree $PT$, a subtree of $PR^{<\omega}$.  An infinite path
  $P$ in $PT$ consists of a set $X$ in which $A$ is recursive and a
  coding of $X'$.  Since $\emptyset'$ can compute a path in $PT$
  uniformly, $\emptyset'$ can uniformly compute a low set $T$-above
  $A$ together with the lowness index for that set.
  
  Trying to keep our presentation simple, we view the function $F$ as
  defined not on $\omega$ but on $PR^{< \omega}$, i.e. on $p$-strings.
  Define
  \[
  F^*(k) = max\{F(\rho) : \rho \in PR^{< \omega} \ \& \ |\rho|
  = k \}.
  \]
  Let $G$ be a total recursive function such that the sequence
  $\{G(x,s)\}_s$  has the limit $F^*(x)$, for all $x$.
  By replacing functions $F^*$ and $G$ with possibly larger functions we may assume the following:
\begin{itemize}
 \item the sequence $\{G(x,s)\}_s$ is nondecreasing
 \item $F^*(x)$ is greater or equal than the modulus
       of this limit, for all $x$, i.e. $j \geq F^*(x)$ implies $G(x,j) = F^*(x)$ for all $j, x$
 \item the function $F^*$ is increasing.
\end{itemize}

  With any $p$-string $\rho = (\sigma,\alpha)$ we will effectively
  associate a recursive tree $Tr_\rho \subseteq 2^{<\omega}$ (see
  below).  We may assume, without loss of generality, that $F^*$ is
  growing sufficiently fast so that if such tree is finite then the
  value of $F^*(|\rho|)$ is greater than the maximal $d$ such that
  $\sigma$ is $d$-extendable on this tree $Tr_\rho$.  That is to say
  that if $Tr_\rho$ is a finite tree, then $F^*(|\rho|)$ is at least
  as large as the height of $Tr_\rho$ above $\sigma$.

  We now build, recursively in $\emptyset'$, an infinite perfect tree
  $PT$ of $p$-strings.  We ensure that for any infinite path $Z =(X,Y)$ on $PT$,
  $X \in 2^\omega, Y \in \{-1,0,1\}^{\omega}$, and if $Z$ is recursive
  in $\emptyset'$, then   
  \[
  \mbox{$X$ is low, $Y^{\geq 0} = X'$ and $A \leq_T X$.}
  \]

  We will build tree $PT$ inductively by stages.  Let $S_0$ consists
  of the empty $p$-string (denoted by $\Lambda$).  At stage $e > 0$,
  we will produce a finite collection of $p$-strings $S_e$ by
  extending $p$-strings from $S_{e-1}$ and we will restrict our tree
  $PT$, at stage $e$, to those $p$-strings compatible with $p$-strings
  from $S_e$.  Thus, we define a sequence of finite trees, ordered
  under end-extension, and let $PT$ be their union.

  At stage $e+1$, we will simultaneously and continuously in $PT$
  decide the $e$-th instance of the jump operator on the first
  coordinates $X$ of the infinite paths $(X,Y)$ through $PT$.

  We fix some notation. Let $DTr_e$, the divergence tree for $\Phi_e$,
  denote the recursive tree
  \[
  DTr_e=\{\sigma \in 2^{<\omega} : \Phi_{e,|\sigma|}(\sigma)(e) \uparrow\}
  \]
  and let $\mathcal D_e = [DTr_e]$ be the set of functions $f$ such
  that $\Phi_e(f)(e)\uparrow$.

  Finally, let $TrPA(A)$ denote an $A$-recursive tree $\subseteq
  2^{<\omega}$ such that $[TrPA(A)] = \mathcal {PA}(A)$.  That is to
  say, the class of infinite paths through $TrPA(A)$ is the class of
  all $\{0,1\}$-valued $A$-DNR functions. Given an infinite recursive
  set $M$, let $CdTr(A,M)$ denote the subtree of $2^{<\omega}$
  recursive in $A$ defined by
  \[
  CdTr(A,M)=\{\sigma \in 2^{<\omega}: Restr(\sigma,M) \in TrPA(A)\}.
  \]
  Thus,
  \[
  [CdTr(A,M)] = \{f :  Restr(f,M) \in \mathcal {PA}(A)\}.
  \]
  Here, $CdTr$ stands for a coding tree.  Observe, that any member of
  $[CdTr(A,M)]$ is T-above $A$.  This is our way of coding $A$ into
  any $X$ such that, for some $Y$, $(X,Y)$ is an infinite path through
  the tree $PT$.

  With each $p$-string $\rho = (\sigma,\alpha)$, we effectively
  associate a recursive tree $Tr_\rho$ $\subseteq 2^{<\omega}$ and a
  $\Pi^0_1$ class $\mathcal B_{\rho} = [Tr_{\rho}]$ in the following
  way.  Let first $Tr_{\Lambda}$ denote a recursive tree such that
  $[Tr_{\Lambda}]$ equals $\mathcal {PA}$.  Further, if $\beta =
  \alpha^{\geq 0}$, then $Tr_{(\sigma,\alpha)} = (Tr_{\Lambda} \cap
  Ext(\sigma)) \bigcap_{\beta(j)=0} DTr_j$.  Intuitively, $\mathcal
  B_{(\sigma,\alpha)} = [Tr_{(\sigma,\alpha)}]$ is a restriction of
  $\mathcal {PA} \cap [\sigma]$ to a $\Pi^0_1$ class of sets $X$ for
  which $j \notin X'$ for any $j$ for which $\beta(j)=0$.  

  We will ensure that for any $p$-string $(\sigma,\alpha) \in S_{e+1}$
  one of the following conditions holds.
  \begin{itemize}
  \item $e \in X'$ for every $X \in [Tr_{(\sigma,\alpha)}]$
  \item $e \notin X'$ for every $X \in [Tr_{(\sigma,\alpha)}]$
  \end{itemize}
  In addition, with each $\rho \in PT$ we associate
  (recursively in $\emptyset'$) an infinite recursive set $M_\rho$.
  Each set $M_\rho$ represents a way of coding of $A$ into (first
  coordinate of) infinite paths extending $\rho$ through $PT$.  We
  will prove that along any infinite path through $PT$
  there will be only finitely many changes in the set that is so associated.
  In other words, along each path $(X,Y)$ our coding of $A$ into $X$
  will stabilize and ensure that $X\geq_T A$.

  In outline, we begin by letting $M_{\Lambda}$ be an infinite recursive
  set such that
  $Restr(\mathcal {PA},M_{\Lambda}) = 2^\omega$.  In
  other words, we commit ourselves to building a 
  $PA$ set (i.e. a $\{0,1\}$-valued DNR function) and we fix an infinite 
  set $M_\Lambda$ for coding $A$. At 
  stage $e+1$ in our construction, if $\rho = (\sigma,\alpha) \in PR
  \cap S_e$ and we can extend $\rho$ during stage $e+1$ without injury
  (as described below), then we will associate the same infinite
  recursive set $M_{\rho}$ with the extensions of $\rho$ that we add
  in $S_{e+1}$.  Otherwise, our construction may be injured at $\rho$.
  In this case, in order to fulfill our commitments about $X$ and $X'$
  which are specified by $\rho$, we must abandon the set $M_{\rho}$ as
  the place to code $A$.  We then specify  new infinite recursive sets
  $M^+$.  The technical device of the construction is to maintain the
  ability to keep numbers out of $X'$ as specified by $\rho$ as we
  monitor the coding of $A$ into the extensions of $\rho$ in $PT$.

  With each $p$-string $\rho$ and infinite recursive set $M$, we
  define the $A$-recursive tree
  \[
  Tr_\rho(A,M) = Tr_{\rho} \cap CdTr(A,M).
  \]
  Similarly, we let
  \[
  \begin{aligned}
    \mathcal {B}_{\rho}(A,M) &= [Tr_{\rho}(A,M)]\\
    &= \mathcal{B}_{\rho} \cap \{f : Restr(f,M) \in \mathcal {PA}(A) \}
  \end{aligned}
  \]
  Intuitively, if $\rho = (\sigma,\alpha)$, then $\mathcal
  {B}_{\rho}(A,M)$ is a restriction of $\mathcal {PA} \cap [\sigma]$ :
  first to the class of sets $X$ for which $j \notin X'$ for any $j$
  for which $\beta(j)=0$ and, second, to the class of sets $X$ for
  which $Restr(X,M) \in \mathcal {PA}(A)$.

  For any $p$-string $(\sigma,\alpha)$ which is on our tree $PT$, we
  will ensure the following two properties.
  \begin{itemize}
  \item $\omega$-extendability of $\sigma$ on $Tr_{(\sigma,\alpha)}$
    (i.e. in $\mathcal B_{(\sigma,\alpha)}$) 
  \item $F^*(|(\sigma,\alpha)|)$-extendability of $\sigma$ on
    $Tr_{(\sigma,\alpha)}(A,M_{(\sigma,\alpha)})$ 
  \end{itemize}
  As already indicated earlier, we suppose that $F^*$ grows
  sufficiently fast so that $Tr_{(\sigma,\alpha)}$ is finite if and
  only if $\sigma$ is not $F^*(|(\sigma,\alpha)|)$-extendable on tree
  $Tr_{(\sigma,\alpha)}$.

 To help the reader we give here a concentrated list of notations most frequently used in what follows.\\
 $PR$ denotes $\{0,1\} \times \{-1,0,1\}$, 
 $PR^{< \omega}$ denotes the set of all finite sequences of elements of $PR$, i.e. the set of $p$-strings.\\
 $PT$ denotes an infinite perfect tree recursive in $\emptyset'$, a subtree of $PR^{< \omega}$, (to be built).\\
 $DTr_e$ denotes the divergence tree for $\Phi_e$, i.e. the recursive tree
 $DTr_e=\{\sigma \in 2^{<\omega} : \Phi_{e,|\sigma|}(\sigma)(e) \uparrow\}$, 
 and $\mathcal D_e = [DTr_e]$ denotes the set of functions $f$ such
 that $\Phi_e(f)(e)\uparrow$.\\
 $TrPA(A)$ denotes an $A$-recursive tree $\subseteq
 2^{<\omega}$ such that $[TrPA(A)] = \mathcal {PA}(A)$, i.e. the class of all 
 infinite paths through $TrPA(A)$ is the class of
 all $\{0,1\}$-valued $A$-DNR functions. \\
 $CdTr(A,M)$, for an infinite recursive set $M$,  ({\it a coding tree}), denotes the subtree of $2^{<\omega}$
 recursive in $A$ defined by $CdTr(A,M)=\{\sigma \in 2^{<\omega}: Restr(\sigma,M) \in TrPA(A)\}$.
 Thus, the class of all infinite paths through this tree is just the class 
 of all functions $f \in 2^\omega$ for which $Restr(f,M) \in \mathcal {PA}(A)$ 
 (i.e. functions coding in a recursive projection given by $M$ a function from $\mathcal {PA}(A)$).\\
 $Tr_\rho$, for each $p$-string $\rho = (\sigma,\alpha)$, denotes a recursive subtree of $2^{<\omega}$
 and $\mathcal B_{\rho}$ denotes a $\Pi^0_1$ class $\mathcal B_{\rho} = [Tr_{\rho}]$ such that 
 $Tr_{(\sigma,\alpha)} = (Tr_{\Lambda} \cap Ext(\sigma)) \bigcap_{\beta(j)=0} DTr_j$,
 where $Tr_{\Lambda}$ is a recursive tree with $[Tr_{\Lambda}] =\mathcal {PA}$
 and $\beta = \alpha^{\geq 0}$. I.e., $\mathcal B_{(\sigma,\alpha)}$ is equal to 
 $\mathcal {PA} \cap [\sigma] \cap \{X : \beta(j)=0 \rightarrow j \notin X' \}.$\\
 For a $p$-string $\rho$ and an infinite recursive set $M$,
 $Tr_\rho(A,M)$ denotes the $A$-recursive tree $Tr_{\rho} \cap CdTr(A,M)$
 and $\mathcal {B}_{\rho}(A,M)$ denotes 
 $[Tr_{\rho}(A,M)]$, i.e. $\mathcal {B}_{\rho}(A,M)$ is equal to 
 $\mathcal{B}_{\rho} \cap \{f : Restr(f,M) \in \mathcal {PA}(A) \}$.
 Intuitively, if $\rho = (\sigma,\alpha)$, then $\mathcal
 {B}_{\rho}(A,M)$ is a restriction of $\mathcal {PA} \cap [\sigma]$ :
 first to the class of sets $X$ for which $j \notin X'$ for any $j$
 for which $\beta(j)=0$ and, second, to the class of sets $X$ for
 which $Restr(X,M) \in \mathcal {PA}(A)$.

 Now, we present the precise recursion step of our construction.

  \subsection{{Stage $e+1$}}

  Let $\rho = (\sigma,\alpha)$ be a $p$-string from $S_e$.  We
  consider several cases.

  \begin{case} At least one of $\sigma*j$, for $j=0,1$, is both
    \begin{itemize}
    \item $\omega$-extendable on $Tr_{\rho} \cap DTr_e$ (i.e. in
      $\mathcal B_{\rho} \cap \mathcal D_e$)
    \item and $F^*(|\rho|+1)$-extendable on $Tr_{\rho}(A,M_{\rho}) \cap
      DTr_e$.
    \end{itemize}
    Then for all such $j$'s put $(\sigma*j,\alpha*0)$ into $S_{e+1}$
    and let $M_{(\sigma*j,\alpha*0)}$ be $M_{\rho}$ .  In this case,
    we have ensured that $e \notin X'$, without injury.
  \end{case}

  \begin{case} 
    The previous case does not apply.
    
    First observe, that necessarily $\mathcal B_{\rho}(A,M_{\rho})
    \cap \mathcal D_e$ is empty (otherwise we would have the previous
    case).  Further, this condition is recursively recognized relative
    to $\emptyset'$.  However, since we are not working with a low
    index of $A$, we cannot determine recursively in
    $\emptyset'$ whether $\mathcal B_{\rho}(A,M_{\rho})$ is empty.  We
    consider two subcases. Either it is possible to ensure $e \in X'$
    without injury at this case or we detect an injury.  Injuries will
    be explained in detail below.

    If we take for each $j=0,1$ all strings
    $\tau, \tau \succeq \sigma*j$, of length $|\rho|+ 1 +
    F^*(|\rho|+1)$, then the only such strings which are on
    $Tr_{\rho}(A,M_{\rho})$ and are $\omega$-extendable in $\mathcal
    B_{\rho}$ (if there are such at all) are not $\omega$-extendable on
    $DTr_e$, (i.e. $[\tau] \cap \mathcal D_e$ is empty).  So, for
    $j=0,1$ take all $\gamma, \gamma \succeq \sigma*j$ (if such exist at
    all), which satisfy the following conditions.

    \begin{itemize}
    \item $\gamma\in Tr_{\rho}(A,M_{\rho})$
    \item $\gamma$ is $\omega$-extendable on $Tr_\rho$ (i.e. extends
      to an element of $\mathcal B_{\rho}$)
  \item $\gamma$ has length $\leq |\rho|+ 1 + F^*(|\rho|+1)$
  \item $[\gamma] \cap \mathcal D_e$ is empty
    \end{itemize}

    Now we split into two subcases, depending on whether there is a
    string $\gamma$ as above which is sufficiently extendable on
    $Tr_{\rho}(A,M_{\rho})$.
    
    \textit{Subcase 2.a.\ \ } There are a string $\gamma$ and $j$,
    such that $0 \leq j \leq 1$, $\gamma \succeq \sigma*j$, $\gamma$
    is both $F^*(|\gamma|)$-extendable on $Tr_{\rho}(A,M_{\rho})$ and
    $\omega$-extendable on $Tr_{\rho}$ (i.e. in $\mathcal B_{\rho}$),
    $[\gamma] \cap \mathcal D_e$ is empty and $|\gamma| \leq |\rho|+ 1
    + F^*(|\rho|+1)$.

    Then for any such $\gamma, j$ for which no $\tau, \sigma*j \preceq
    \tau \prec \gamma$ has this property, put $(\gamma,\alpha*(-1)^k*1)$
    into $S_{e+1}$, where $k = |\gamma| - |\sigma*j| \ (= |\gamma| -
    |\rho| - 1)$, and let also $M_{\tau} = M_{\rho}$ for any $\tau, \rho
    \prec \tau \preceq (\gamma,\alpha*(-1)^k*1)$.

    In this case, we have ensured that $e \in X'$, without injury.
    
    {\it Subcase 2.b}.  Now assume that neither of the two previous
    situations applies.  Then, we are unable to respect our
    commitments to deciding the jump while continuing the coding of
    $A$.  This is the case in which we injure our coding strategy.

    Observe, that $\mathcal B_{\rho}(A,M_{\rho})$ must be empty.  This
    means that $\sigma$ is only finitely-extendable on
    $Tr_{\rho}(A,M_{\rho})$.  Further, from $\rho$ and using $A$, we
    can compute an upper bound of this finite-extendability.  Thus,
    there are strings $\gamma$, $\gamma \succeq \sigma$, with the following properties.
    \begin{itemize}
    \item $\gamma$ is $\omega$-extendable on $Tr_{\rho}$ (i.e. in
      $\mathcal B_{\rho}$)
    \item $F^*(|\gamma|)$-extendable on $Tr_{\rho}(A,M_{\rho})$
    \end{itemize}
    but, however, neither of the immediate extensions of $\gamma$ is both
    $\omega$-extendable on $Tr_{\rho}$ (i.e. in $\mathcal B_{\rho}$)
    and $F^*(|\gamma|+1)$-extendable on $Tr_{\rho}(A,M_{\rho})$.  
    There are only finitely many of these strings $\gamma$.  In
    particular, each is less than or equal to $|\rho|+ 1 +
    F^*(|\rho|+1)$.  Finally, note that $Tr_{\rho}(A,M_{\rho})$ is a
    subtree of $Tr_{\rho}$.

    Then for each such $\gamma$, we say that an injury occurred at
    $(\gamma,\alpha*(-1)^{k})$, where $k = |\gamma| - |\sigma| \ (=
    |\gamma| - |\rho|)$, and we do the following.

    Let for $j = 0,1, \ d_j$ denote the maximal $d$ such that
    $\gamma*j$ is $d$-extendable on $Tr_{\rho}(A,M_{\rho})$.  Since
    $d_j < F^*(|\gamma|+1)$, let $t_0$ be the least $t$ such that
    $G(|\gamma|+1,t) > d_j$ for both $j=0,1$.  Intuitively, for
    strings $\gamma$ in the current situation, at step
    $t_0$ the recursive approximation of $F^*$ by $G$ is able to see
    that an injury occurred.

    We also know that the maximal $d$ for which $\gamma$ is
    $d$-extendable on $Tr_{\rho}(A,M_{\rho})$ is greater than
    $F^*(|\gamma|)$.  Denote it as $d_\emptyset$.  We now position
    ourselves to use the hypothesis that $F$ eventually dominates
    every function which is partial recursive relative to $A$.  We will
    use $d_\emptyset$ to define a value of an $A$-partial recursive
    function at input $(\gamma,\alpha*(-1)^k)$, which is greater than
    the corresponding value of $F^*$.  Precisely, the value of the
    defined function will be greater than the value $F^*(|\gamma|)$
    and, therefore, also greater than $F(\gamma,\alpha*(-1)^k)$, since
    $F^*(|\rho|) \geq F(\rho)$ for all $p$-strings $\rho$.

    To summarize, $d_\emptyset > F^*(|\gamma|)$, but $d_j <
    F^*(|\gamma|+1)$, $j=0,1$ and also $d_\emptyset = 1 +
    max(d_0,d_1)$.  Further, there is at least one string $\tau, \tau
    \succ \gamma$, which is $\omega$-extendable on
    $Tr_{(\gamma,\alpha*(-1)^k)}$ (i.e. in $\mathcal
    B_{(\gamma,\alpha*(-1)^k)}$ ) and for which $|\tau| - |\gamma| =
    t_0$.  Our ongoing commitments concerning the jumps of the paths
    in $PT$ can be enforced on the extensions of these strings $\tau$.

    For each such $\tau$, we have one of the following two possibilities.
    \begin{enumerate}
    \item The first possibility is that $[\tau] \cap \mathcal D_e =
      \emptyset$.  Then for $q = k + t_0 = |\tau| - |\sigma|$ denote
      $(\tau,\alpha*(-1)^{q-1}*1)$ by $\xi$, put $\xi$ into $S_{e+1}$,
      and let $M_\eta = M_\rho$ for any $p$-string $\eta, \rho \prec
      \eta \prec \xi$.  Finally, apply Lemma~\ref{2.6} to effectively find
      an infinite recursive set $M^+$ such that $Restr(\mathcal
      B_\xi,M^+) = 2^\omega$ and let $M_\xi = M^+$.
    \item The second possibility is that $[\tau] \cap \mathcal D_e
      \neq \emptyset$.  Then for $q = k + t_0 = |\tau| - |\sigma|$
      denote $(\tau,\alpha*(-1)^{q-1}*0)$ by $\xi$, put $\xi$ into
      $S_{e+1}$, and let $M_\eta = M_\rho$ for any $p$-string $\eta, \rho
      \prec \eta \prec \xi$.  Again, apply Lemma~\ref{2.6} to
      effectively find an infinite recursive set $M^+$ such that
      $Restr(\mathcal B_\xi,M^+) = 2^\omega$ and let $M_\xi = M^+$.
    \end{enumerate}
    In both of these possibilities, we start with a new version of
    coding of $A$ into (first coordinate of) infinite paths through
    $PT$ extending $\xi$.
  \end{case}

 This ends the action of our construction during stage $e+1$.
  
 Recall that $PT$ is built inductively by stages.
 At stage $e > 0$ we have produced a finite collection of $p$-strings $S_e$ by extending 
 $p$-strings from $S_{e-1}$ and we restricted our tree
 $PT$ at this stage to those $p$-strings compatible with $p$-strings
 from $S_e$. It is important to note that at any stage $e > 0$
 each $p$-string from $S_{e-1}$ is really extended to at least one (possibly more)
 $p$-string from $S_e$. This fact follows immediately from our construction. 
 Thus, $PT$ is an infinite perfect tree.

  \subsection{Verification}

  It remains to show that our construction achieves its aims.  We must
  show that for any infinite path $(X,Y)$ through $PT$ which is
  computable from $\emptyset'$, $Y^{\geq 0} = X'$ and $A \leq_T X$.
  Since $(X,Y)$ is recursive relative to $\emptyset'$ and $X'$ is
  recursive in $Y$, $X$ is low as required.

  It is clear that during stage $e+1$ we have decided the membership
  of $e$ in $X'$ for any such $(X,Y)$ and so $Y^{\geq 0}=X'$.  More
  precisely, using $\emptyset'$ we can find $\rho = (\sigma,\alpha)
  \in S_{e+1}$, $\sigma \prec X$ and then $e \in X'$ if and only if
  $\alpha(|\rho|) = \alpha(|\alpha|) = 1$.  It remains only to verify
  $A \leq_T X$.  For that it is sufficient to show that along any
  infinite path through $PT$ there are only finitely many stages where
  an injury occurs.

  For this purpose, we build a partial function $H$ on $PR^{<\omega}$
  (i.e. on $p$-strings) recursively in $A$.  In the definition of $H$,
  we $A$-recursively approximate the $\emptyset'$-construction of
  $PT$.  At the beginning, all $p$-strings are associated with an
  infinite recursive set $M_\Lambda$ (the original set used for coding
  into $\mathcal {PA}$).  During building $H$ we sometimes either stop
  some strategies or restart strategies for defining values of $H$ by
  changing current values of infinite recursive sets associated with
  $p$-strings, more precisely, $p$-strings preceding $(\sigma,\alpha)$
  may eventually either stop a strategy for defining
  $H(\sigma,\alpha)$ or restart a strategy for defining
  $H(\sigma,\alpha)$ (if such strategy was not already finished) by
  changing a current value of $M$ associated with $(\sigma,\alpha)$.

  We will now first describe an isolated strategy of $H$ for a
  $p$-string $(\sigma,\alpha)$ and an infinite recursive set $M$ and
  then we will indicate how to combine strategies together.

  Consider a tree $Tr_{(\sigma,\alpha)}(A,M)$.  If $\sigma$ is
  $\omega$-extendable on $Tr_{(\sigma,\alpha)}(A,M)$ then the strategy
  has no output and no effect on $p$-strings extending
  $(\sigma,\alpha)$.  If $\sigma$ is not $\omega$-extendable on
  $Tr_{(\sigma,\alpha)}(A,M)$, then define $H(\sigma,\alpha) =
  d_\emptyset$, where $d_\emptyset$ denotes the maximal $d$ such that
  $\sigma$ is $d$-extendable on $Tr_{(\sigma,\alpha)}(A,M)$. Further,
  wait for a step $t$ such that $G(|\sigma|+1,t) > d_j$ for both
  $j=0,1$, where $d_j$ denotes the maximal $d$ for which $\sigma*j$ is
  $d$-extendable on $Tr_{(\sigma,\alpha)}(A,M)$ (recall that
  $d_\emptyset = 1 + max(d_0,d_1)$).  If there is no such $t$, the
  strategy has no effect on $p$-strings extending $(\sigma,\alpha)$.
  If there is such $t$, take the least such and denote it by $t_0$.
  Take a finite collection $Q$ of all $p$-strings
  $(\tau,\alpha*(-1)^{(t_0-1)}*j)$, for $j=0,1$ and $\tau \succ
  \sigma$ with $|\tau| - |\sigma| = t_0$, and stop all strategies
  (which are still active) for $p$-strings extending $(\sigma,\alpha)$
  either with length $< |\sigma|+t_0$, or not compatible with any
  $p$-string from $Q$, and finally restart strategies for any
  $p$-string $\eta$ of length $\geq |\sigma| + t_0$ compatible with
  some $p$-string from $Q$ but now with a new infinite recursive set
  $M$ which is determined as follows. For any $p$-string $\xi =
  (\tau,\alpha*(-1)^{(t_0-1)}*j)$ from $Q$, as we did in our
  construction, effectively find an infinite recursive set $M^+$ such
  that if $\mathcal B_\xi \neq \emptyset$ then $Restr(\mathcal
  B_\xi,M^+) = 2^\omega$, and $\xi$ together with all $p$-strings
  extending such $\xi$ are restarted with this newly associated set
  $M^+$.

  All strategies for defining values of $H$ are combined together
  easily by a finite injury style where any $p$-string has a higher
  priority than exactly all $p$-strings extending it.  As mentioned
  previously, at the beginning all $p$-strings are associated with
  $M_\Lambda$ (the original set used for coding into $\mathcal {PA}$).

  We omit further details.

  By virtue of its definition, $H$ is an $A$-partial recursive
  function. The only values of $H$ that are relevant for our purposes
  are those on $p$-strings from tree $PT$.
  
  It remains to verify that $A \leq_T X$ for any infinite path $(X,Y)$
  through $PT$.  Let such path $(X,Y)$ be given.  It is clearly
  sufficient to show that there are only finitely many $p$-strings
  $\rho_i, \rho_i \prec (X,Y)$, at which an injury occurs. Suppose,
  for a contradiction, that $\{\rho_i\}_i$ is an infinite sequence of
  $p$-strings with increasing length at which an injury occurred and
  for which $\rho_i \prec (X,Y)$.  It follows from the construction
  and our assumptions on $F^*$ and $G$ that $H(\rho_i)$ is defined and
  greater than $F^*(|\rho_i|)$ for all $i$.  Since $F^*(|\rho|) \geq
  F(\rho)$ for any $p$-string $\rho$, it immediately yields a failure
  of $F$ to eventually dominate all $A$-partial recursive functions. A
  contradiction.
  
%

 To finish the construction we have to provide an infinite path $(X,Y)$ through $PT$
 which is recursive in $\emptyset'$, and then let $A^* = X$. This is easy. Recall that $PT$ 
 is an infinite perfect tree which is recursive in $\emptyset'$. Thus, e.g. we may fix 
 the left-most (in the standard ordering) infinite path through $PT$. This path is obviously 
 recursive in $\emptyset'$.

 We end with a remark.  The above construction can be
 carried out uniformly in a $\emptyset'$-index of $A$.  That is,
 there are recursive functions $f,g$ with properties stated in Lemma~{3.5}.

\end{proof}

\begin{remark}
  \begin{enumerate}
  \item In the proof of Lemma~\ref{3.5} we have coded a set $A$ into
    $A^*$ not straightforwardly.  We could do that, but instead, we
    have nested into a $\Pi^0_1$ class $\mathcal {PA}$ a $\Pi^{0,A}_1$
    class $\mathcal {PA}(A)$ all members of which are $T$-above
    $A$. This way does not increase the complexity of the proof and is
    more general. It can be used to provide just one global
    construction of a $T$-upper bound mentioned in Theorem~\ref{3.2}.
\item It is easy to verify that for any infinite path $(X,Y)$ through
  $PT$ (not necessarily computable by $\emptyset'$) in the above proof
  $A \leq_T X$, $X' = Y^{\geq 0}$ and $X \in GL_1$.
  \end{enumerate}
\end{remark}

Before giving a proof of Theorem~\ref{3.2} we
illustrate the main idea (of a more difficult implication) on a simple example.

\begin{claim}\label{3.6} 
  Given a sequence $\{A_n\}_n $ of low sets uniformly recursively in
  $\emptyset'$ so that low indices of all finite joins of members of
  the sequence are also uniformly recursive in $\emptyset'$ then there
  is a low set $A$ which is $T$-above all $A_n$.  Moreover, we may
  require $\bf a >> \bf a_n$ for all n, where $\bf a, \bf a_n$ are
  $T$-degrees of $A$, $A_n$ respectively.
\end{claim}

\begin{proof}
  We give a sketch of the main idea.  Let $B_n = \bigoplus_{i=0}^{n-1}
  A_i$ for $n \geq 1$ \ and $B_0 = \emptyset$.  We construct
  recursively in $\emptyset'$ a sequence of nonempty relativized
  $\Pi^0_1$ classes $\{\mathcal B_n\}_n$ such that $\mathcal B_n
  \supseteq \mathcal B_{n+1}$ for all $n$ and the required low set $A$
  recursive in $\emptyset'$ will be the only set in the intersection
  of all these classes.  Each class $\mathcal B_n$ is a nonempty
  $\Pi^{0,B_n}_1$ class.  We also construct recursively in
  $\emptyset'$ a sequence of infinite recursive sets $\{M_e\}_e$ such
  that for all $e$, $M_e \supseteq M_{e+1}$ and $Restr(\mathcal
  B_e,M_e) = \mathcal {PA}(B_e)$.

  Let $\mathcal B_0$ be $\mathcal {PA}$ and $M_0 = \omega$. So, since $B_0 = \emptyset$, 
  $Restr(\mathcal B_0,M_0) = \mathcal {PA}(B_0)$.

  At step $e+1$ we first use the relativized Low Basis Theorem of
  \citet{Jockusch.Soare:1972} to force the $e$-th instance of the jump
  operator.  This gives a subclass $\mathcal B^*_{e+1} \subseteq
  \mathcal B_e$.  Obviously $Restr(\mathcal B^*_{e+1},M_e) \subseteq
  \mathcal {PA}(B_e)$ and $Restr(\mathcal B^*_{e+1},M_e)$ is also a
  nonempty $\Pi^{0,B_e}_1$ class.  Thus, by a modification of
  Lemma~\ref{2.6}, there is an infinite recursive set $M_{e+1}$, a
  subset of $M_e$, such that $Restr(\mathcal B^*_{e+1},M_{e+1}) =
  2^\omega$.  Let $\mathcal B_{e+1} = \mathcal B^*_{e+1} \cap \{f:
  Restr(f,M_{e+1}) \in \mathcal {PA}(B_{e+1})\}$.  Then, $\mathcal
  B_{e+1}$ is a $\Pi^{0,B_{e+1}}_1$ class and $Restr(\mathcal
  B_{e+1},M_{e+1}) = \mathcal {PA}(B_{e+1})$.

  This ends step $e+1$.

  One verifies that there is just only one set $A$ in the intersection
  of all these classes and that $A$ satisfies all required
  requirements.  Moreover, $T$-degree of $A$ is \ $>>$ \ $T$-degrees
  of all $A_n$.
\end{proof}

\begin{remark} We could equally construct, recursively in $\emptyset'$, an infinite perfect tree
$(\subseteq 2^{<\omega} )$ such that any infinite path through it both 
is $T$-above all $A_n$ and belongs to $GL_1$. We could also work with elements of
$2^\omega \times 2^\omega$ and build, recursively in $\emptyset'$, an infinite perfect tree 
$(\subseteq (\{0,1\} \times \{0,1\})^{<\omega})$ such that for any infinite path $(X,Y)$ 
through it $Y$ codes (in fact, is equal to) the jump of $X$.
\end{remark}

However, in case we do not have low indices of such low sets $A_n$ available 
(recursively in $\emptyset'$) and we have only a weaker information about lowness of these sets,  
we have to adapt the method explained in the proof of Claim~\ref{3.6} by combining it 
with the technique of the proof of Lemma~\ref{3.5}.
We explain subsequently how to do it.

\begin{lemma}\label{3.7} Given a low set $B$ with a low index $b$ and a function $F$ 
  recursive in $\emptyset'$, there is a uniform
  way to obtain from a $\emptyset'$-index of a set $A$ with the
  property that any partial function recursive in $B \oplus A$ is eventually dominated by
  $F$ both a low set $A^*$ and an index of lowness of $A^*$ such that
  $B \oplus A \leq_T A^*$.  That is to say that there are recursive functions
  $f, g$ such that if $\Phi_e(\emptyset')$ is total and equal to some
  set $A$ so that any partial function recursive in $B \oplus A$ is eventually dominated
  by $F$ then $\Phi_{f(e)}(\emptyset')$ is a low set, $g(e)$ is its
  lowness index and $B \oplus A \leq_T \Phi_{f(e)}(\emptyset')$.
\end{lemma}

\begin{proof}
It is an easy relativization of Lemma~\ref{3.5}, i.e. the construction in the proof
of Lemma~\ref{3.5} is relativized to an oracle $B$ and it begins now with the class
$\mathcal {PA}(B)$  instead of the class $\mathcal {PA}$. 
Let us explicitely note that we substantially use a low index $b$ of $B$ so that 
we can always decide whether a $B$-recursive tree $(\subseteq 2^{<\omega})$ is finite 
or infinite. Thus, we can always keep $\omega$-extendability of relevant strings in 
infinite $B$-recursive trees, similarly as we kept $\omega$-extendability of relevant strings 
in infinite recursive trees in the original construction (recall that, in the proof of Lemma~\ref{3.5}, \ 
$(\sigma, \alpha) \in PT$ implies $\omega$-extendability of $\sigma$ in $Tr_{(\sigma, \alpha)}$).
\end{proof}

The situation changes (slightly) if instead of a low index of a low set $B$ we have only a weaker information 
about lowness of $B$, namely, if we have a function recursive in $\emptyset'$ which eventually dominates all partial function 
recursive in $B$.

\begin{lemma}\label{3.8} Given a function $F$ recursive in $\emptyset'$ and a low set $B$ there is a uniform
  way to obtain from a $\emptyset'$-index of a set $A$ with the
  property that any partial function recursive in $B \oplus A$ is eventually dominated by
  $F$ both a low set $A^*$ and an index of lowness of $A^*$ such that
  $B \oplus A \leq_T A^*$.  That is to say that there are recursive functions
  $f, g$ such that if $\Phi_e(\emptyset')$ is total and equal to some
  set $A$ so that any partial function recursive in $B \oplus A$ is eventually dominated
  by $F$ then $\Phi_{f(e)}(\emptyset')$ is a low set, $g(e)$ is its
  lowness index and $B \oplus A \leq_T \Phi_{f(e)}(\emptyset')$.
\end{lemma}

\begin{proof}
  It is a further generalization of the technique used in the proof of
  Lemma~\ref{3.5} and Lemma~\ref{3.7}.  Suppose that any partial
  function recursive in $B \oplus A$ is eventually dominated by $F$.
  We have to replace missing low index of $B$ by a technique of
  Lemma~\ref{3.5}, i.e. we use guesses computed by $F$ (or $F^*$) as
  approximations to correct answers to $\Sigma^{0,B}_1$ questions.  By
  a finite injury style, these approximations are eventually really
  correct answers to such questions from some point on. From such a
  point the construction resembles that one described in the proof of
  Lemma~\ref{3.7}. It means, the construction then simulates the proof
  of Lemma~\ref{3.7} with a deeper nested class $\mathcal {PA}(A)$
  into the current version of $\mathcal {PA}(B)$.  Thus, to combine
  everything together, we nest relativized $\Pi^0_1$ classes in a
  cascade style (like in Claim~\ref{3.6}), here with depth of nesting
  equal to $2$.  We first nest a $\Pi^{0,B}_1$ class $\mathcal
  {PA}(B)$ into $\mathcal {PA}$ and then we further nest a $\Pi^{0,A}$
  class $\mathcal {PA}(A)$ into it (by applying Lemma~\ref{2.6} like
  in the proof of Claim~\ref{3.6}).
  %
  %
  Injury occurs at a stage $s$ to one of the $\Pi^0_1$ guesses carried
  forward from the previous stage when the construction encounters a
  witness to the effect that the corresponding $\Pi^0_1$-sentence is
  false.  At such a stage, the construction restarts all guesses on
  trees which are nested at or below the first such injured one.
  These constitute the first guess which requires attention and those
  of lower priority.  The construction does not change the remaining,
  higher-priority, guesses.  The crucial thing is that at the highest
  level (level $0$), where we work with nonrelativized $\Pi^0_1$ class
  $\mathcal {PA}$, we never have any injury because $\emptyset'$ can
  always correctly decide any $\Sigma^0_1$ question.  Under the given
  assumptions the construction reaches all required goals.  The
  verification is a standard application of the finite injury
  technique and we omit further details.
\end{proof}

\begin{proof}
  We now complete the proof of Theorem~\ref{3.2}.  As it
  was already mentioned, the implication from (2) to (1) is direct.
    
  To verify the implication from (1) to (2), we may identify the ideal
  $\mathcal C$ with the ideal $\mathcal A$ described in (1).  
  Given a function $F$ recursive in $\emptyset'$ and a sequence of sets $\{A_n\}_n$
  which satisfy given assumptions, we combine the method of the proof of Claim~\ref{3.6} 
  with the method of the proof of Lemma~\ref{3.8}, i.e. missing low indices are replaced by the function $F$ (or $F^*$).
  At each step we have to nest a next relativized $\Pi^0_1$ class into a previous one, where nested classes are 
  $\mathcal {PA}(A_n)$ \ (or $\mathcal {PA}( \bigoplus^{n-1}_{j=0} A_j)$). Thus, the depth 
  of nesting of relativized $\Pi^0_1$ classes increases by $1$ at each step. Analogously as before, 
  at each level of nesting there are only finitely many injuries. We omit further details.

\end{proof}

\section{A question}

We have left the following questions open.

\begin{question}
  Is there a natural set of conditions which characterize whether an
  ideal in the Turing degrees of the recursively enumerable sets has
  a upper bound which is low and recursively enumerable?
\end{question}

 Exact pairs play an important role in the study of
 $T$-degree structures. By a result of \citet{nerode.shore:1980} it
 follows that there is an exact pair for the class of $K$-trivials in
 $\Delta^0_2$ $T$-degrees.  On the other hand the following problem is
 open.

 \begin{question}
   Is there a low exact pair for the class of $K$-trivial sets,
   i.e. are there low sets $A,B$ such that $\mathcal K = \{C : C \leq_T
   A \ \& \ C \leq_T B \}$? 
 \end{question}

 Finally, since the $K$-trivial sets are closely linked with the
 1-random reals, the following special situation is interesting.

 \begin{question}
   Is there a low $1$-random set which is a $T$-upper bound for the
   class of $K$-trivial sets?
 \end{question}


\begin{thebibliography}{25}
\providecommand{\natexlab}[1]{#1}
\providecommand{\url}[1]{\texttt{#1}}
\expandafter\ifx\csname urlstyle\endcsname\relax
  \providecommand{\doi}[1]{doi: #1}\else
  \providecommand{\doi}{doi: \begingroup \urlstyle{rm}\Url}\fi

\bibitem[Ambos-Spies and Ku{\v{c}}era(2000)]{Ambos-Spies.Kucera:2000}
K.~Ambos-Spies and A.~Ku{\v{c}}era.
\newblock Randomness in computability theory.
\newblock In \emph{Computability theory and its applications (Boulder, CO,
  1999)}, volume 257 of \emph{Contemp. Math.}, pages 1--14. Amer. Math. Soc.,
  Providence, RI, 2000.

\bibitem[Chaitin(1977)]{Chaitin:1977}
G.~J. Chaitin.
\newblock Algorithmic information theory.
\newblock \emph{IBM J. Res. Develop}, 21:\penalty0 350--359, 496, 1977.

\bibitem[Downey and Hirschfeldt(2007)]{Downey.Hirschfeldt:nd}
R.~Downey and D.~R. Hirschfeldt.
\newblock Algorithmic randomness and complexity.
\newblock to appear, 2007.

\bibitem[Downey et~al.(2006)Downey, Hirschfeldt, Nies, and
  Terwijn]{Downey.Hirschfeldt.ea:2006}
R.~Downey, D.~R. Hirschfeldt, A.~Nies, and S.~A. Terwijn.
\newblock Calibrating randomness.
\newblock \emph{Bull. Symbolic Logic}, 12\penalty0 (3):\penalty0 411--491,
  2006.
\newblock ISSN 1079-8986.

\bibitem[Downey et~al.(2003)Downey, Hirschfeldt, Nies, and
  Stephan]{Downey.Hirschfeldt.ea:2003}
R.~G. Downey, D.~R. Hirschfeldt, A.~Nies, and F.~Stephan.
\newblock Trivial reals.
\newblock In \emph{Proceedings of the 7th and 8th Asian Logic Conferences},
  pages 103--131, Singapore, 2003. Singapore Univ. Press.

\bibitem[Hirschfeldt et~al.(2007)Hirschfeldt, Nies, and
  Stephan]{Hirschfeldt.Nies.ea:nd}
D.~R. Hirschfeldt, A.~Nies, and F.~Stephan.
\newblock Using random sets as oracles.
\newblock to appear, 2007.

\bibitem[Jockusch and Soare(1972)]{Jockusch.Soare:1972}
C.~G. Jockusch, Jr. and R.~I. Soare.
\newblock {$\Pi_1^0$} classes and degrees of theories.
\newblock \emph{Trans. Amer. Math. Soc.}, 173:\penalty0 33--56, 1972.


\bibitem[Ku{\v{c}}era(1989)]{Kucera:1989}
A.~Ku{\v{c}}era.
\newblock On the use of diagonally nonrecursive functions.
\newblock In \emph{Logic Colloquium '87 (Granada, 1987)}, volume 129 of
  \emph{Stud. Logic Found. Math.}, pages 219--239. North-Holland, Amsterdam,
  1989.

\bibitem[Ku{\v{c}}era(1993)]{Kucera:1993}
A.~Ku{\v{c}}era.
\newblock On relative randomness.
\newblock \emph{Ann. Pure Appl. Logic}, 63\penalty0 (1):\penalty0 61--67, 1993.
\newblock ISSN 0168-0072.
\newblock 9th International Congress of Logic, Methodology and Philosophy of
  Science (Uppsala, 1991).

\bibitem[Ku{\v{c}}era and Terwijn(1999)]{Kucera.Terwijn:1999}
A.~Ku{\v{c}}era and S.~A. Terwijn.
\newblock Lowness for the class of random sets.
\newblock \emph{J. Symbolic Logic}, 64\penalty0 (4):\penalty0 1396--1402, 1999.
\newblock ISSN 0022-4812.


\bibitem[Li and Vit{\'a}nyi(1997)]{Li.Vitanyi:1997}
M.~Li and P.~Vit{\'a}nyi.
\newblock \emph{An introduction to {K}olmogorov complexity and its
  applications}.
\newblock Graduate Texts in Computer Science. Springer-Verlag, New York, second
  edition, 1997.
\newblock ISBN 0-387-94868-6.

\bibitem[Miller and Nies(2006)]{Miller.Nies:2006}
J.~S. Miller and A.~Nies.
\newblock Randomness and computability: open questions.
\newblock \emph{Bull. Symbolic Logic}, 12\penalty0 (3):\penalty0 390--410,
  2006.
\newblock ISSN 1079-8986.

\bibitem[Nerode and Shore(1980)]{nerode.shore:1980}
A.~Nerode and R.~A. Shore.
\newblock Reducibility orderings: theories, definability and automorphisms.
\newblock \emph{Ann. Math. Logic}, 18:\penalty0 61--89, 1980.

\bibitem[Nies(2005)]{Nies:2005}
A.~Nies.
\newblock Lowness properties and randomness.
\newblock \emph{Adv. Math.}, 197\penalty0 (1):\penalty0 274--305, 2005.
\newblock ISSN 0001-8708.

\bibitem[Nies(2006)]{Nies:2006}
A.~Nies.
\newblock Reals which compute little.
\newblock In \emph{Logic Colloquium '02}, volume~27 of \emph{Lect. Notes Log.},
  pages 261--275. Assoc. Symbol. Logic, La Jolla, CA, 2006.

\bibitem[Nies(2007)]{Nies:nd*d}
A.~Nies.
\newblock Computability and randomness.
\newblock To appear, 2007.

\bibitem[Odifreddi(1989)]{Odifreddi:1989}
P.~Odifreddi.
\newblock \emph{Classical recursion theory}, volume 125 of \emph{Studies in
  Logic and the Foundations of Mathematics}.
\newblock North-Holland Publishing Co., Amsterdam, 1989.
\newblock ISBN 0-444-87295-7.
\newblock The theory of functions and sets of natural numbers, With a foreword
  by G. E. Sacks.

\bibitem[Odifreddi(1999)]{Odifreddi:1999}
P.~Odifreddi.
\newblock \emph{Classical recursion theory. {V}ol. {II}}, volume 143 of
  \emph{Studies in Logic and the Foundations of Mathematics}.
\newblock North-Holland Publishing Co., Amsterdam, 1999.
\newblock ISBN 0-444-50205-X.

\bibitem[Robinson(1971)]{Robinson:1971}
R.~W. Robinson.
\newblock Interpolation and embedding in the recursively enumerable degrees.
\newblock \emph{Ann. of Math.}, 93:\penalty0 285--314, 1971.

\bibitem[Scott(1962)]{Scott:1962}
D.~Scott.
\newblock Algebras of sets binumerable in complete extensions of arithmetic.
\newblock In \emph{Recursive Function Theory}, volume~5 of \emph{Proceedings of
  Symposia in Pure Mathematics}, pages 117--121, Providence, R.I., 1962.
  American Mathematical Society.

\bibitem[Simpson(1977)]{Simpson:1977}
S.~Simpson.
\newblock Degrees of unsolvability: a survey of results.
\newblock In J.~Barwise, editor, \emph{Handbook of mathematical logic}, pages
  631–--652. North-Holland, Amsterdam, 1977.

\bibitem[Soare(1987)]{Soare:1987}
R.~I. Soare.
\newblock \emph{Recursively Enumerable Sets and Degrees}.
\newblock Perspectives in Mathematical Logic, Omega Series. Springer--Verlag,
  Heidelberg, 1987.

\bibitem[Yates(1969)]{Yates:1969}
C.~E.~M. Yates.
\newblock On the degrees of index sets. {II}.
\newblock \emph{Trans. Amer. Math. Soc.}, 135:\penalty0 249--266, 1969.
\newblock ISSN 0002-9947.

\end{thebibliography}

\end{document}